\documentclass[a4paper,oneside,12pt]{article}
\usepackage{amssymb,amsmath,amsfonts,amsthm}
\usepackage{graphicx,graphics}
\usepackage{epsfig}
\usepackage{latexsym}
\usepackage[applemac]{inputenc}
\usepackage{ae,aecompl}
\usepackage{mathrsfs}
\usepackage[english]{babel}
 \usepackage[colorlinks=true]{hyperref}
\usepackage{pstricks}
\newcommand{\E}[1]{\mathbf{E} \left [ #1 \right ]}
\newcommand{\sun}{\mathbb{S}_{1}}

\renewcommand{\P}[1]{\mathbf{P}\left( #1 \right)}

\renewcommand{\geq}{\geqslant}
\renewcommand{\leq}{\leqslant}

\newcommand{\op}[1]{\operatorname{#1 }}
 \usepackage[colorlinks=true]{hyperref}

\usepackage[english]{babel}
\newtheorem{theorem}{Theorem}[section]
\newtheorem*{theorem*}{Theorem}

\newtheorem{proposition}[theorem]{Proposition}
\newtheorem{lemma}[theorem]{Lemma}

\newtheorem{rek}[theorem]{Remark}

\title{Random laminations and multitype branching processes}

\author{Nicolas Curien\footnote{Dma-Ens, 45 rue d'Ulm 75005 Paris, France (e-mail: nicolas.curien@ens.fr)} \  and Yuval Peres\footnote{Microsoft Research, Redmond, Washington, USA (e-mail: peres@microsoft.com)}}
\begin{document}

\maketitle
\abstract{We consider multitype branching processes arising in the study of random laminations of the disk. We classify these processes according to their subcritical or supercritical behavior and provide Kolmogorov-type estimates in the critical case corresponding to the random recursive lamination process of \cite{CLG10}. The proofs use the infinite dimensional Perron-Frobenius theory and quasi-stationary distributions.}
\section{Introduction} 
\begin{figure}[!h]
\begin{center}
\includegraphics[width=16cm]{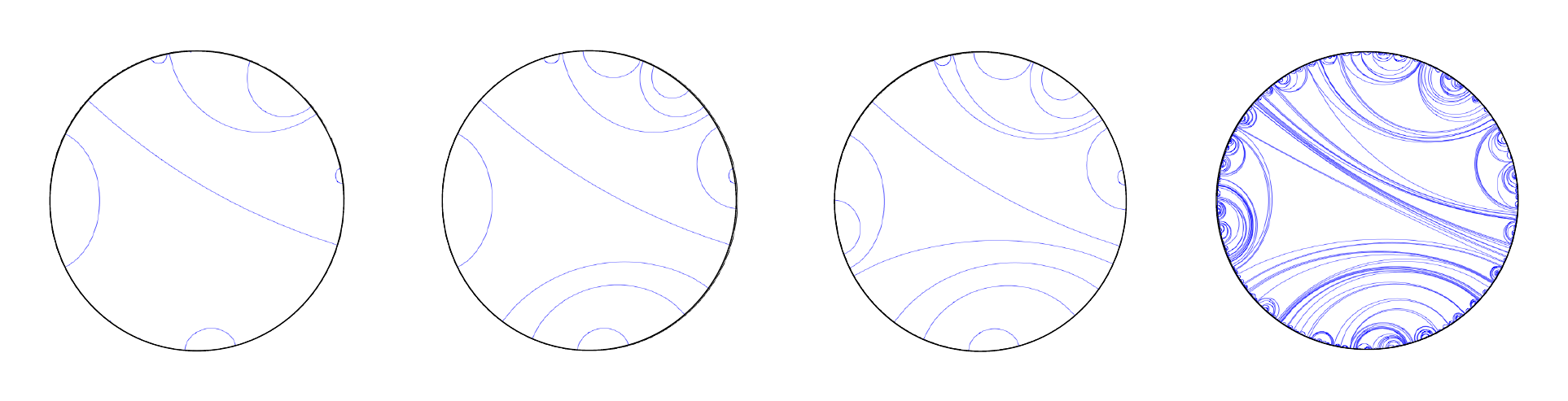}
\end{center}
\caption{\label{brw:fig1}An illustration of the process creating the sequence $(L_n)_{n\geq1}$. We use hyperbolic chords rather than Euclidean chords for aesthetic reasons.}

\end{figure}

In this note we are interested in multitype branching processes that arise in the study of random recursive laminations. In order to introduce and motivate our results, let us briefly recall  the basic construction of \cite{CLG10}. Consider a sequence $U_1,V_1,U_2,V_2,\ldots$ of independent
random variables, which are uniformly distributed over the unit circle $\sun$. We then construct inductively
a sequence $L_1,L_2,\ldots$ of random closed subsets of the closed unit disk $\overline{\mathbb D}$.
To start with, $L_1$ is set to be the (Euclidean) chord $[U_1V_1]$ with endpoints $U_1$ and $V_1$. Then at step $n+1$, we consider two cases. Either the chord $[U_{n+1}V_{n+1}]$
intersects $L_n$, and we put $L_{n+1}=L_n$. Or the chord $[U_{n+1}V_{n+1}]$ does not
intersect $L_n$, and we put $L_{n+1}=L_n \cup [U_{n+1}V_{n+1}]$. Thus, for every integer $n\geq 1$, 
$L_n$ is a disjoint union of random chords. See Fig.\,\ref{brw:fig1}.

\begin{figure}[!h]

\begin{center}
\includegraphics[width=13cm]{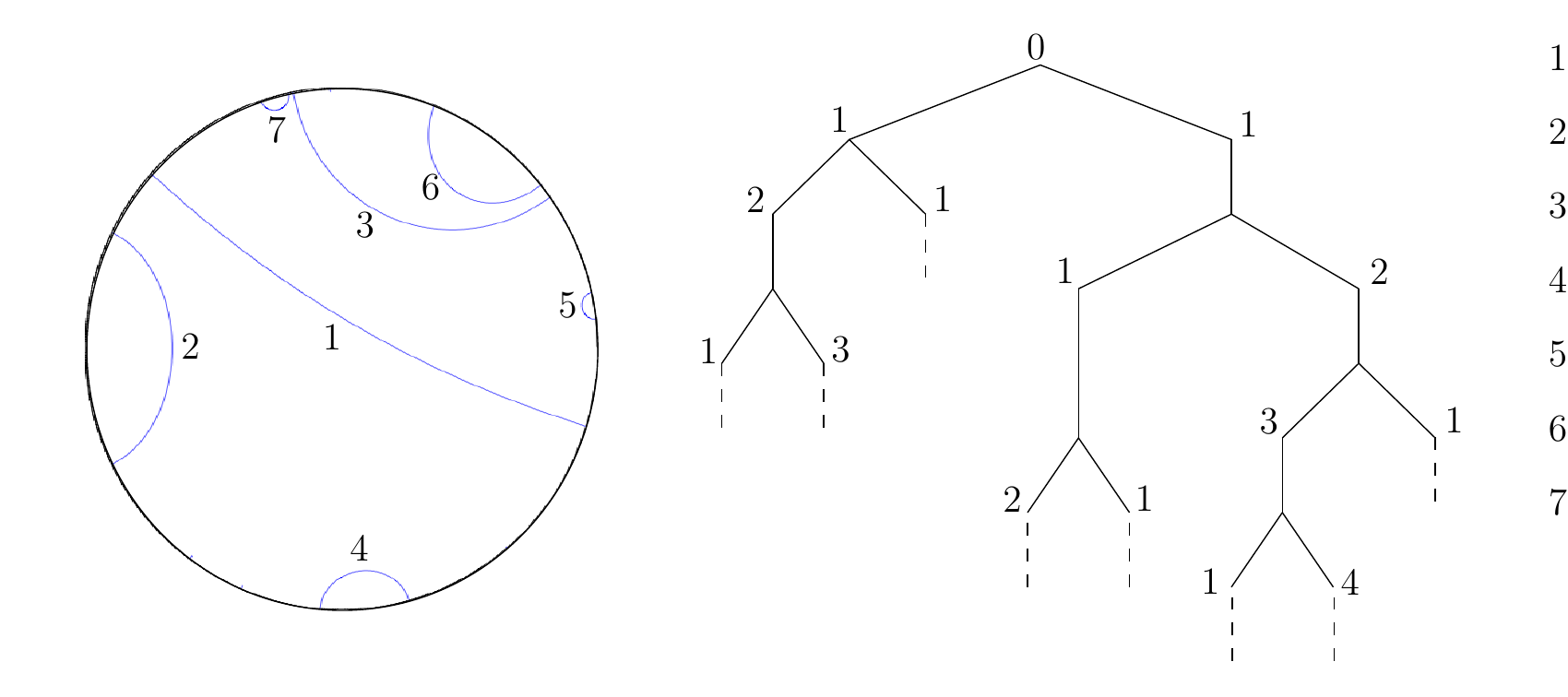}
\end{center}
\caption{\label{brw:fig2}On the left-hand side, the first $7$ chords of the splitting process. On the right-hand side, the associated branching process corresponding to the number of ends of the fragments at their creations. Notice that we split the fragments according to the order of appearance of the chords, thus the binary tree on the right-hand side seems stretched.}
\end{figure}

 A fragment of $L_{n}$ is a connected component of $\overline{\mathbb{D}} \backslash L_{n}$.  These fragments have a natural genealogy that we now describe. The first fragment, $\overline{\mathbb{D}}$, is represented by $\varnothing$. Then the first chord $[U_{1}V_{1}]$ splits $\overline{\mathbb{D}}$ into two fragments, which are viewed as the offspring of $\varnothing$. We then order these fragments in a random way: With probability $1/2$, the first child of $\varnothing$, which is represented by $0$, corresponds to the largest fragment and the second child, which is represented by $1$, corresponds to the other fragment. With probability $1/2$ we do the contrary. We then iterate this device (see Fig.\,\ref{brw:fig2}) so that each fragment appearing during the splitting process is labeled by an element of 
the infinite binary tree
$$\mathbb{T}_{2}=\bigcup_{n \geq 0} \{0,1\}^n\,,\ \quad \mbox{ where } \{0,1\}^0= \{\varnothing\}. $$

If $F$ is a fragment, we call {\it end} of $F$, any connected component of
$F\cap \sun$. 
For convenience, the full disk $\overline{\mathbb{D}}$ is viewed as a fragment with $0$ end. Consequently, we can associate to any $u \in \mathbb{T}_{2}$ a label $\ell(u)$ that corresponds to the number of ends of the corresponding fragment in the above process. Lemma 5.5 of \cite{CLG10} then entails that this random labeling of $\mathbb{T}_{2}$ is described by the following branching mechanism: For any $u \in \mathbb{T}_{2}$ labeled $m \geq 0$, choose $m_1 \in \{0,1, \ldots , m \}$ uniformly at random and assign the values $1+m_1$ and $1+m-m_1$ to the two children of $u$. This is the multitype branching process we will be interested in. See Fig.\,\ref{brw:fig2}.

We can also define a random labeling by using the above branching mechanism but starting with a value $a \geq 0$ at the root $\varnothing$ of $\mathbb{T}_{2}$, the probability distribution of this process will be denoted $\mathbf{P}_{a}$ and its relative expectation $\mathbf{E}_{a}$. A \emph{ray} is an infinite geodesic path $\mathbf{u}= (u_1,u_2, \ldots) \in \{0,1\}^{\mathbb{N}}$ starting from the root $\varnothing$ in $\mathbb{T}_2$. For any ray $\mathbf{u} =(u_1, \ldots , u_n,\ldots)$ or any word of finite length $u=(u_1, \ldots , u_n)$, we denote  by $[\mathbf{u}]_i$ or $[u]_i$   the word $(u_1, \ldots , u_i)$ for $1\leq i \leq n$, and $[u]_{0}=\varnothing$.

 \begin{theorem*}[{\cite[Lemma 5.5]{CLG10}}] \label{before} Almost surely, there exists no ray $\mathbf{u}$ along which all the labels starting from $4$ are bigger than or equal to $4$, 
\begin{eqnarray*}\mathbf{P}_{4}\big(\exists \mathbf{u} \in \{0,1\}^\mathbb{N} : \ell([\mathbf{u}]_i) \geq 4, \forall i \geq 0 \big) &=& 0.\end{eqnarray*}
\end{theorem*}
The starting label $4$ does not play any special role and can be replaced by any value bigger than $4$. This theorem was proved and used in \cite{CLG10} to study certain properties of the random closed subset $L_{\infty}= \overline{\cup L_{n}}$, and in particular to prove that it is  almost surely a \emph{maximal} lamination (roughly speaking that the complement of $L_\infty$ is made of disjoint triangles), see \cite[Proposition 5.4]{CLG10}. One of the purposes of this note is to provide quantitative estimates related to this theorem. Specifically let $$G_n = \big\{ u \in \{0,1\}^n : \ell([u]_i) \geq 4, \forall i \in \{0,1, \ldots , n\}\big\}$$ be the set of paths in $\mathbb{T}_2$ joining the root to the level $n$ along which the labels are bigger than or equal to $4$.

\begin{theorem} \label{main} The expected number of paths starting from the root and reaching level $n$ along which the labels starting from $4$ are bigger than or equal to $4$ satisfies 
\begin{eqnarray} \mathbf{E}_{4}\left[\#G_n\right] &\underset{n \to \infty}{\longrightarrow}& \frac{4}{e^2-1}. \label{first} \end{eqnarray} Furthermore, there exist two constants $0<c_1<c_2< \infty$ such that the probability that $G_n \ne \varnothing$ satisfies
\begin{eqnarray} \frac{c_1}{n} \quad \leq & \mathbf{P}_{4}\big(G_n \ne \varnothing \big) &\leq \quad  \frac{c_2}{n}.\end{eqnarray}
\end{theorem}
\begin{rek} These estimates are reminiscent of the critical case for Galton-Watson processes with finite variance $\sigma^2 < \infty$. Indeed if $H_n$ denotes the number of vertices at height $n$ in such a process then $\E{H_n}=1$ and Kolmogorov's estimate \cite{Kol38} implies that $\P{H_n\ne 0} \sim \frac{2}{\sigma^2 n}.$
\end{rek}

\begin{figure}[!h]
\begin{center}
\includegraphics[width=17cm]{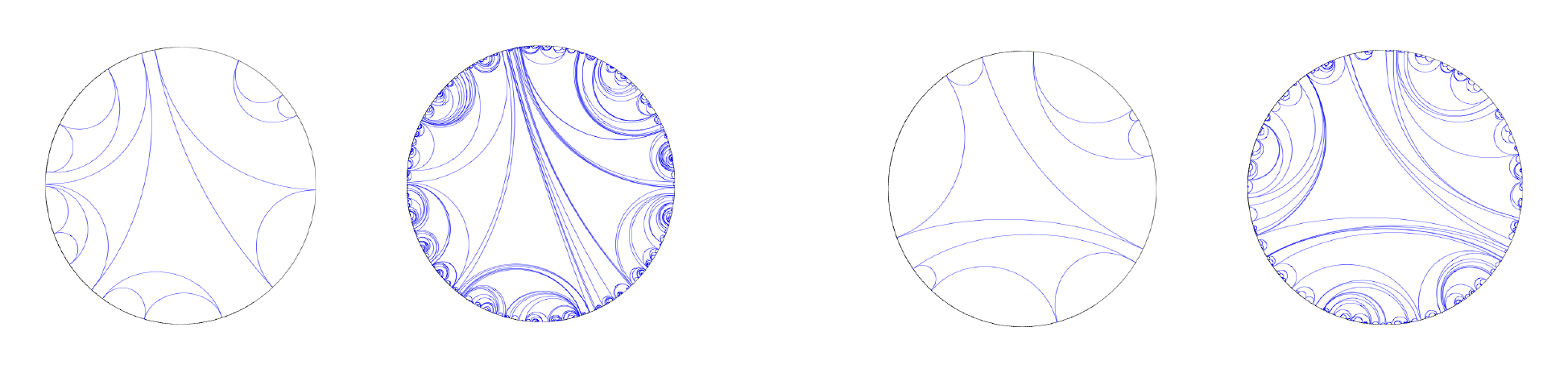}
\end{center}
\caption{\label{brw:fig3}Extension of the process $(L_n)_{n\geq1}$ where we throw triangles or squares instead of chords.}
\end{figure}

The proof of Theorem \ref{main} relies on identifying the quasi-stationary distribution of the labels along a fixed ray conditioned to stay bigger than or equal to $4$. This is done in Section \ref{sec:critical}. In Section $3$, we also study analogues of this branching random walk on the $k$-ary tree, for $k\geq 3$, coming from a natural generalization of the process $(L_n)_{n \geq 0}$ where we replace chords by triangles, squares... see Fig.\,\ref{brw:fig3}. 

We prove in these cases that there is no critical value playing the role of $4$ in the binary case.

 \smallskip
 {\bf Acknowledgments.} The first author thanks Microsoft Research and the University of Washington, where most of this work was done, for their hospitality. We are also grateful to Jean-Fran\c cois Le Gall for precious comments and suggestions on a first version of this note.

\section{The critical case}
\label{sec:critical}
\subsection{A martingale} \label{martingale}
Fix an arbitrary ray $\mathbf{u}_0$ in $\mathbb{T}_2$, for example $\mathbf{u}_0= (0,0,0,0,0, \ldots)$ and define $X_n= \ell([\mathbf{u}_0]_n)$ for $n \geq 0$, so that $X_{n}$ is the value at the $n$-th vertex on the fixed ray $\mathbf{u}_0$ of the $\mathbb{T}_2$-indexed walk $\ell$ starting from $x_0 \geq 4$ at the root. Then $(X_n)_{n \geq 0}$ is a homogeneous Markov chain with transition probabilities given by 
$$ P_{2}(x,y) = \frac{1}{x+1} \mathbf{1}_{1 \leq y \leq x+1}.$$ We first recall some results derived in \cite{CLG10}. If $\mathcal{F}_n$ is the canonical filtration of $(X_n)_{n \geq 0}$ then a straightforward calculation leads to $\mathbf{E}_{x_{0}}[X_{n+1}\mid \mathcal{F}_n] = 1 + X_n/2$, hence the process $M_n= 2^n(X_n-2)$ is a martingale starting from $x_0-2$. For $i \geq 1$, we let $T_{i}$ be the stopping time $T_{i}=\inf\{n \geq 0: X_{n}=i\},$ and $T= T_{1}\wedge T_{2}\wedge T_{3}$. By the stopping theorem applied to the martingale $(M_{n})_{n\geq0}$, we obtain for every $n \geq 0$,
$$ x_0-2 =\mathbf{E}_{x_{0}}[M_{n \wedge T}] =\mathbf{E}_{x_{0}}[-2^{T_{1}}\mathbf{1}_{\{T_{1}=T \leq n\}}] +\,0\,+  
\mathbf{E}_{x_{0}}[2^{T_{3}}\mathbf{1}_{\{T_{3}=T \leq n\}}] +\mathbf{E}_{x_{0}}[2^n(X_{n}-2)\mathbf{1}_{\{T>n\}}].$$ One can easily check from the transition kernel of the Markov chain $(X_{n})_{n\geq 0}$ that for every $i \geq 1$, 
$ \mathbf{P}_{x_{0}}[T_{1}=T=i]= \mathbf{P}_{x_{0}}[T_{2}=T=i]=\mathbf{P}_{x_{0}}[T_{3}=T=i]$
. Hence, the equality in the last display becomes 
$$ x_0-2= \mathbf{E}_{x_{0}}[2^n(X_{n}-2)\mathbf{1}_{\{T>n\}}],$$ or equivalently
\begin{eqnarray}x_0-2 &=& 2^n \mathbf{P}_{x_{0}}[T>n]\; \mathbf{E}_{x_{0}}[X_{n}-2\mid T>n]. \label{condprob} \end{eqnarray}
Our strategy here is to compute the stationary distribution of $X_n$ conditionally on the non extinction event $\{T>n\}$, in order to prove the convergence of $\mathbf{E}_{4}[X_n\mid T>n]$ and finally to get asymptotics for $\mathbf{P}_{4}[T>n]$. Before any calculation, we make a couple of simple remarks. Obviously $\mathbf{E}_{x_{0}}[X_{n}-2 \mid T>n]\geq 2$, and thus we get $2^n\,\mathbf{P}_{x_{0}}(T>n) \leq \frac{x_0-2}{2}$. Since there are exactly $2^n$ paths joining the root $\varnothing$ of $\mathbb{T}_2$ to the level $n$, we deduce that the number $\#G_n$ of paths joining $\varnothing$ to the level $n$ along which the labels are bigger than or equal to $4$ satisfies \begin{eqnarray} \mathbf{E}_{x_{0}}[\#G_n] &\leq& \frac{x_0-2}{2}. \label{upper}\end{eqnarray}
Notice that a simple argument shows that if $4\leq x_0\leq x_1$ then the chain $X_n$ starting from $x_0$ and the chain
$X_n'$ starting from $x_1$ can be coupled in such a way that $X_n \leq X_n'$ for all $n \geq 0$.

\subsection{The quasi-stationary distribution}
We consider the substochastic matrix of the Markov chain $X_n$ killed when it reaches $1,2$ or $3$:  This is the matrix $(\tilde{P}_{2}(x,y))_{x,y\geq 4}$ given by
$$ \tilde{P}_{2}(x,y) = \frac{1}{x+1} \mathbf{1}_{y \leq x+1}.$$

We will show that $\tilde{P}_{2}$ is a $2$-recurrent positive matrix, in the sense of \cite[Lemma 1]{SVJ66}. For that purpose we seek left and right non-negative eigenvectors of $\tilde{P}_{2}$ for the eigenvalue $1/2$. In other words we look for two sequences $(g(x))_{x\geq 4}$ and $(f(x))_{x\geq 4}$ of non-negative real numbers such that $f(4)=g(4)=1$ (normalization) and for every $x \geq 4$ 
\begin{eqnarray} g(x) \quad=\quad 2 \sum_{y\geq 4} g(y) \tilde{P}_{2}(y,x) \quad=\quad  2\sum_{ y = (x-1)\vee 4}^{\infty} \frac{g(y)}{y+1}, \label{left} \\
 f(x) \quad=\quad 2\sum_{y \geq 4}\tilde{P}_{2}(x,y) f(y) \quad=\quad\frac{2}{x+1} \sum_{y=4}^{x+1} f(y). \label{right}
\end{eqnarray}
We start with the left eigenvector $g$. From \eqref{left}, we get $g(5)=g(4)=0$,  and $g(i)-g(i+1)=\frac{2}{i}g(i-1)$ for $i \geq 5$. Letting $$ G(z) = \sum_{i\geq 4} \frac{z^{i+1}}{i+1}g(i), \quad  0 \leq z <1,$$
 the last observations lead to the following differential equation for $G$ \begin{eqnarray*} 2G(z) &=& \displaystyle z^{-1}(z-1)G'(z)  + z^3,\end{eqnarray*} with the condition $G(z) = z^5/5+o(z^5)$. A simple computation yields $ G(z) = 3/4 \exp(2z)(z-1)^2 + ({z^3}/{2}+ {3z^2}/{4}- {3}/{4}).$ After normalization, the generating function $G_{1/2}(z) =\sum_{i\geq 4} g_{1/2}(i) z^i$ of the unique probability distribution $g_{1/2}$ which is a left eigenvector for the eigenvalue $1/2$ is given by 
$$G_{1/2}(z)=\frac{z}{2} \big(\exp(2z)(z-1) + z + 1\big),$$ that is 
$$ g_{1/2}(i) = \frac{2^{i-3}(i-3)}{(i-1)!}\mathbf{1}_{i \geq 4}.$$
This left eigenvector is called the \emph{quasi-stationary distribution} of $X_n$ conditioned on non-extinction. For the right eigenvector $f,$ a similar approach  using generating functions is possible, but it is also easy to check by induction that  
 
 $$f(i) = \frac{i-2}{2} \mathbf{1}_{i \geq 4},$$ satisfies \eqref{right}. Hence the condition $(iii)$ of Lemma $1$ in \cite{SVJ66} is fulfilled  and the substochastic matrix $\tilde{P}_{2}$ is $2$-recurrent positive. For every $x\geq 4$, set $q_{n}(x)= \mathbf{P}_{4}(X_{n}=x \mid T>n) = \mathbf{P}_{4}(T>n)^{-1} \tilde{P}_{2}^n v\,(x)$ where $v$ stands for the ``vector'' $(v_{i})_{i\geq4}$ with $v_{4}=1$ and $v_{i}=0$ if $i\geq 5$. Theorem $3.1$ of \cite{SVJ66} then implies that 
  \begin{eqnarray} q_n(x) &\underset{n \to \infty}{\longrightarrow}& g_{1/2}(x). \label{cfinite}\end{eqnarray}
 Unfortunately this convergence does not immediately imply that $\mathbf{E}_{4}[X_n\mid T>n] \longrightarrow \E{X}$ where $X$ is distributed according to $g_{1/2}$. But this will  follow from the next proposition.
 
\begin{proposition} \label{decroissant}For every $n\geq 0$ the sequence 
$ \displaystyle \left( \frac{q_n(x)}{g_{1/2}(x)}\right)_{x \geq 4}$ is decreasing.\end{proposition}
\proof By induction on $n \geq 0$. For $n=0$ the statement is true. Suppose it holds for $n \geq 0$. By the definition of $q_{n+1}$, for $x\geq 4$ we have 
 \begin{eqnarray}
 q_{n+1}(x)&=& \mathbf{P}_{4}(X_{n+1}=x \mid T>n+1)\nonumber \\
  &=& \frac{1}{\mathbf{P}_{4}(T>n+1)}\sum_{z\geq 4} \mathbf{P}_{4}(X_{n}=z \,,\ X_{n+1}=x\,,\ T>n) \nonumber \\
  & =& \frac{\mathbf{P}_{4}(T>n)}{\mathbf{P}_{4}(T>n+1)}\sum_{z \geq (x-1)\vee 4} \frac{q_{n}(z)}{z+1} \label{etoile}\end{eqnarray}
 We need to verify that, for every $x \geq 4$,  we have $q_{n+1}(x)g_{1/2}(x+1) \geq q_{n+1}(x+1)g_{1/2}(x)$ or equivalently, using \eqref{etoile} and \eqref{left} with $g=g_{1/2}$, that
\begin{eqnarray}
\left(\sum_{z \geq x\vee 4} \frac{g_{1/2}(z)}{z+1}\right) \left(\sum_{z \geq (x-1)\vee 4} \frac{q_n(z)}{z+1}\right) & \geq &  \left(\sum_{z \geq (x-1)\vee 4} \frac{g_{1/2}(z)}{z+1}\right) \left(\sum_{z \geq x\vee 4} \frac{q_n(z)}{z+1}\right) \nonumber 
\end{eqnarray}
For $x=4$ this inequality holds. Otherwise, if $x >4$, we have to prove that
\begin{eqnarray}
 q_n(x-1)\sum_{z \geq x \vee 4} \frac{g_{1/2}(z)}{z+1} & \geq & g_{1/2}(x-1)\sum_{z\geq x \vee 4} \frac{q_{n}(z)}{z+1}. \label{etoile2}
\end{eqnarray}
Set $A_{x} = \frac{q_{n}(x-1)}{g_{1/2}(x-1)}$ to simplify notation. The induction hypothesis guarantees that $q_{n}(z) \leq A_{x} g_{1/2}(z)$ for every $z \geq x$, and therefore $$\sum_{z\geq x \vee 4} \frac{q_{n}(z)}{z+1} \quad \leq \quad A_{x} \sum_{z \geq x\vee 4} \frac{g_{1/2}(z)}{z+1}.$$
This gives the bound \eqref{etoile2} and completes the proof of the proposition. \endproof
By Proposition \ref{decroissant} we have for every $x\geq 1$, $ \frac{q_{n}(x)}{g_{1/2}(x)} \leq \frac{q_{n}(4)}{g_{1/2}(4)} \leq C,$ where $C=\sup_{n\geq 0}\frac{q_{n}(4)}{g_{1/2}(4)} <\infty$ by \eqref{cfinite}. This allows us to apply dominated convergence to get $$\mathbf{E}_{4}[X_{n}|T>n] = \sum_{x \geq 4}x q_{n}(x) \xrightarrow[n\to\infty]{} \sum_{x\geq 4}^{}x g_{1/2}(x) =  G_{1/2}'(1) = \frac{e^2+3}{2}.$$ Using \eqref{condprob} we then conclude that 
\begin{eqnarray} 2^n\mathbf{P}_{4}[T>n] &\xrightarrow[n\to\infty]{}& \frac{4}{e^2-1} \label{proba}.
\end{eqnarray}
\subsection{Proof of Theorem \ref{main}}
We first introduce some notation. 
We denote  the tree $\mathbb{T}_2$ truncated at level $n$ by $\mathbb{T}_2^{(n)}$. For every $u=(u_1,\ldots,u_n)\in\{0,1\}^n$, and every $j\in\{0,1,\ldots,n\}$, recall that $[u]_j=(u_1,\ldots,u_j)$, and if $j\geq 1$, also set $[u]_j^*=(u_1,\ldots,u_{j-1},1-u_{j})$. We say that $j \in \{0,1, \ldots , n-1\}$ is a \emph{left turn} (resp.\,\emph{right turn}) of $u$ if $u_{j+1}=0$ (resp. $u_{j+1}=1$). A \emph{down step} of $u$ is a time $j \in \{0,1, \ldots , n-1\}$ such that 
$$ \ell([u]_j)-\ell([u]_{j+1}) \geq 2.$$
Note that if $j$ is a down step of $u$ then $\ell([u]_{j+1}^*) = 2 + \ell([u]_j)-\ell([u]_{j+1}) \geq 4$.  The set of all $j \in \{0,1, \ldots , n-1\}$ that are left turns, resp.\,right turns, resp.\,down steps, of $u$ is denoted by $\op{L}(u)$, resp.\,$\op{R}(u)$, resp\,$\op{D}(u)$. We endow $\mathbb{T}_2$ with the lexicographical order $\preceq$\,,\ and say that a path $u \in \{0,1\}^n$ is on the left (resp.\,right) of $v \in \{0,1\}^n$ if $u \preceq v$ (resp.\,$v\preceq u$).  A vertex of $\{0,1\}^n$ will be identified with the path it defines in $\mathbb{T}_{2}^{(n)}$. If $u,v \in \mathbb{T}_2$ we let $u\wedge v$ be the last common ancestor of $u$ and $v$.

\proof[Proof of Theorem \ref{main}.] \textsc{Lower bound.} We use a second moment method. Recall that $$G_n = \big\{ u \in \{0,1\}^n : \ell([u]_i) \geq 4, \forall i \in \{0,1, \ldots , n\}\big\}$$ is the set of all paths in $\mathbb{T}_2^{(n)}$ from the root to the level $n$ along which the labels are bigger than or equal to $4$. A path in $G_n$ is called "good". Using \eqref{proba}, we can compute the expected number of good paths and get \begin{eqnarray} \nonumber\mathbf{E}_{4}[\#G_n]&=&2^n\mathbf{P}_{4}[T>n] \quad \underset{n\to \infty}{\longrightarrow} \quad \frac{4}{e^2-1},\end{eqnarray} as $n \to \infty$, which proves the convergence \eqref{first} in the theorem. For $u\in G_{n}$ and $j \in \{0,1,\ldots,n\}$, we let $\op{Right}(u,j)$ be the set of all good paths to the right of $u$ that diverge from $u$ at level $j$, 
$$ \op{Right}(u,j) = \{ v \in G_n: u \preceq v \mbox{ and }u \wedge v = [u]_j\}.$$ In particular, if $j$ is a right turn for $u$, that is $u_{j+1}=1$, then $\op{Right}(u,j)=\varnothing$. Furthermore $\op{Right}(u,n)=\{u\}$. Let us fix a path $u \in \{0,1\}^n$, and condition on $u \in G_n$ and on the labels along $u$. Let $j \in \{0,1,2, \ldots,n\}.$ Note that the first vertex of a path in $\op{Right}(u,j)$ that is not an ancestor of $u$ is $[u]_{j+1}^*$ and its label is $2+\ell([u]_{j})-\ell([u]_{j+1})$, so if we want $\op{Right}(u,j)$ to be non-empty, the time $j$ must be a down step of $u$.  If $j$ is a left turn and a down step for $u$, the subtree $\{ w \in \mathbb{T}_2^{(n)}: w \wedge [u]_j^* = [u]_j^*\}$ on the right of $[u]_j$ is a copy of $\mathbb{T}_2^{(n-j-1)}$, whose labeling starts at $\ell([u]_{j+1}^*)$. Hence thanks to \eqref{upper} we get 
\begin{eqnarray*}\mathbf{E}_{4}[\# \op{Right}(u,j) \mid u \in G_n\,,\ (\ell([u]_i))_{0 \leq i\leq n}] &\leq& \frac{\ell([u]_{j+1}^*)-2}{2}= \frac{\ell([u]_{j})-\ell([u]_{j+1})}{2}.\end{eqnarray*}
Since the labels along the ancestral line of $u$ cannot increase by more that one at each step, if $u \in G_{n}$ we have $\sum_{i=0}^{n-1} \mid \ell([u]_{i+1})-\ell([u]_i)\mid \mathbf{1}_{i\in \op{D}(u)} \leq n$. Combining these inequalities, we obtain 
\begin{eqnarray*} \mathbf{E}_{4}\left[\sum_{j=0}^n \# \op{Right}(u,j)  \ \Big| \   u \in G_n\,,\  (\ell([u]_i))_{0 \leq i\leq n}\right] &\leq& \frac{n}{2}.\end{eqnarray*}
We can now bound $\mathbf{E}_{4}[\#G_n^2]$ from above:
\begin{eqnarray} \mathbf{E}_{4}[\#G_n^2] & \leq& 2 \mathbf{E}_{4}\left[\sum_{u \in \{0,1\}^n} \sum_{u \preceq v} \mathbf{1}_{u \in G_n}\mathbf{1}_{v \in G_n}\right]  \nonumber \\
&=& 2\sum_{u \in \{0,1\}^n} \mathbf{P}_{4}(u \in G_n)\mathbf{E}_{4}\left[\sum_{u \preceq v} \mathbf{1}_{v \in G_n} \ \Big|\  u \in G_n\right]  \nonumber \\
& = & 2 \sum_{u \in \{0,1\}^n} \mathbf{P}_{4}(u \in G_n) \mathbf{E}_{4}\left[\sum_{j=0}^n\# \op{Right}(u,j)\ \Big| \  u \in G_n\right]\nonumber \\
&\leq&  n.\label{second} \end{eqnarray} The lower bound of Theorem \ref{main} directly follows from the second moment method : Using \eqref{first} and \eqref{second} we get the existence of $c_{1}>0$ such that 
\begin{eqnarray}  \P{\#G_n >0} \geq \frac{\E{\# G_n}^2}{\E{(\# G_n)^2}} \geq \frac{c_1}{n}. \label{lower}\end{eqnarray}
\textsc{Upper Bound.} We will first provide estimates on the number of down steps of a fixed path $u \in \{0,1\}^n$. Recall that $\op{L}(u)$, $\op{R}(u)$ and $\op{D}(u)$ respectively denote the left turns, right turns, and down steps times of $u$.
\begin{lemma} There exists a constant $c_3>0$ such that, for every $n\geq0$ and every $u_{0} \in \{0,1\}^n$ $$\mathbf{P}_{4}( u_0  \in G_n\,,\ \# \op{D}(u_0) \leq c_3n ) \leq c_3^{-1}2^{-n}\exp(-c_3n).$$
\end{lemma}
\proof We use the notation of Section \ref{martingale}. For any set $A \subset \{ 0,1, \ldots , n-1\}$ and $m \in \{0,1, \ldots , n - \#A\}$, with the notation $N^A_n=\#\{j\in\{0,1,\ldots,n-1\}\backslash A: X_j =5\}$ we have from \cite[formula (27)]{CLG10}
$$ {\mathbf P}\left[X_{j+1}\geq (X_j-1)\vee 4\,,\ \forall j\in\{0,1,\ldots,n-1\}\backslash A \, ,\ N^A_n=m\right]
\leq \Big(\frac{1}{2}\Big)^m
\Big(\frac{3}{7}\Big)^{n-m-\#A},$$
 We will first obtain crude estimates for $N^A_n$. Note that $N^
A_n \leq N^\varnothing_n$ and that $\sup_{i \geq 1}P_2(i,5) = \frac{1}{5}$, so that for any $B \subset \{0,1, \ldots , n\}$ we have 
$$ \P{X_i = 5\,,\ \forall i \in B} \leq 5^{-\#B}.$$ By summing this bound over all choices of $B$ with $\#B\geq m$ we get  $\P{N_n^\varnothing \geq m} \leq 2^n 5^{-m}$ for every $m \in \{0,1, \ldots n\}$. Let $\kappa_{1}\in (0,1/2)$ and $\kappa_{2} \in (0,1)$ such that $\kappa_{1}+\kappa_{2}<1$. We have 
\begin{eqnarray}  &&\P{u_0  \in G_n \,,\ \# \op{D}(u_0) \leq \kappa_{1}n}\nonumber \\  & \leq & \P{u_0 \in G_n\,,\ \# \op{D}(u_0) \leq \kappa_{1}n\,,\   N_n^\varnothing \leq \kappa_{2}n} + \P{N_n^\varnothing \geq \kappa_{2}n} \nonumber \\
& \leq & \sum_{\begin{subarray}{c} A \subset \{0,1, \ldots , n-1\} \\  \#A \leq \kappa_{1}n \end{subarray}} 
{\mathbf P}\left[X_{j+1}\geq (X_j-1)\vee 4,\ \forall j\in\{0,1,\ldots,n-1\}\backslash A\; ; N^A_n\leq \kappa_{2}n\right]
 + \P{N_n^\varnothing \geq \kappa_{2}n} \nonumber \\
& \leq & (\lfloor \kappa_{2}n \rfloor+1)\sum_{\begin{subarray}{c} A \subset \{0,1, \ldots , n-1\} \\
 \#A \leq \kappa_{1}n \end{subarray}} \left(\frac{7}{6}\right)^{\lfloor\kappa_{2}n \rfloor} \left(\frac{3}{7}\right)^{n-\lfloor\kappa_{1}n\rfloor} + \P{N_n^\varnothing \geq \kappa_{2}n} \nonumber \\
 &\leq & n {n \choose \lfloor \kappa_{1}n \rfloor}\left(\frac{7}{6}\right)^{\lfloor\kappa_{2}n\rfloor} \left(\frac{3}{7}\right) ^{n-\lfloor\kappa_{1}n\rfloor} + 2^n 5^{-\lfloor\kappa_{2}n \rfloor} \label{etoile3}
 \end{eqnarray}
Notice that for every $A>1$ we can choose $\kappa_{1}>0$ small enough so that $n {n \choose \lfloor \kappa_{1}n \rfloor} \leq A^n$ for $n$ large enough. Furthermore $$\left(\frac{7}{6}\right)^{\lfloor\kappa_{2}n\rfloor} \left(\frac{3}{7}\right) ^{n-\lfloor\kappa_{1}n\rfloor}=2^{-n}2^{\lfloor \kappa_{1}n\rfloor} \left(\frac{6}{7}\right)^{n-\lfloor \kappa_{1}n \rfloor - \lfloor \kappa_{2}n \rfloor},$$ and by choosing $\kappa_{1}$ even smaller if necessary we can ensure that the right hand side of \eqref{etoile3}  is bounded by $c_{3}^{-1}2^{-n} \exp(-c_{3}n)$ for some $c_{3}>0$.
\endproof
We use the last lemma to deduce that 
\begin{eqnarray}n \mathbf{P}_{4}\left( \exists u \in G_n\,,\  \# \op{D}(u) \leq c_3n\right) &\leq& \frac{n}{c_3}\exp(-c_3n) \underset{n \to \infty}{\longrightarrow} 0. \label{troppetit}\end{eqnarray}
We now argue on the event $E_L=\{ \exists u \in G_n\,,\  \# (\op{D}(u) \cap \op{L}(u)) \geq c_3n/2\}$. On this event there exists a path $u \in G_n$ with at least  $c_3n/2$ down steps which are also left turns. Conditionally on this event we consider  the left-most path $P$ of $G_n$ satisfying these properties, that is 
$$ P = \min_{\preceq} \big\{ u \in G_n\,,\  \# (\op{D}(u) \cap \op{L}(u)) \geq c_3n/2\big\}.$$
A moment's thought shows that conditionally on $P$ and on the values of the labels along the ancestral line of $P$,  the subtrees of $\mathbb{T}_2^{(n)}$ hanging on the right-hand side of $P$, that are the offsprings of the points $[P]_{j+1}^*$ for $j \in \op{L}(P)$,  are independent and distributed as labeled trees started at $\ell([P]_{j+1}^*)$.  

Hence conditionally on $P$ and on the labels $((\ell([P]_{i}), 0 \leq i \leq n)$, for any $j \in \op{L}(P) \cap \op{D}(P)$ the expected number of paths belonging to the set $\op{Right}(P,j)$ (defined in the proof of the lower bound) is \begin{eqnarray} \mathbf{E}_{4}\left[\#\op{Right}(P,j) \ \Big| \ P\,,\  (\ell([P]_i))_{0 \leq i \leq n}\right] &=& 2^{n-j-1}\mathbf{P}_{\ell([P]_{j+1}^*)}\left(T>n-j-1\right) \nonumber \\ &\geq& 2^{n-j-1}\mathbf{P}_{4}(T>n-j-1)\nonumber \\ &\geq& \kappa_{3}>0,\end{eqnarray} where $\kappa_{3}$ is a positive constant independent of $n$ whose existence follows from \eqref{proba}. Thus we have 
\begin{eqnarray}
\mathbf{E}_{4}[\#G_n \mid E_L] &=&\mathbf{E}_{4}\left[\E{\sum_{j=0}^n \# \op{Right}(P,j) \ \Big | \  P\,,\ (\ell([P]_i))_{0 \leq i \leq n} } \ \Big |\  E_L\right] \nonumber \\
& \geq &\kappa_{3}\mathbf{E}_{4}[\# (\op{D}(P) \cap \op{L}(P)) \mid E_L] . \nonumber \\
& \geq& \frac{c_3\kappa_{3}}{2}n. \end{eqnarray}
Since $\mathbf{P}_{4}(E_L) \leq \mathbf{E}_{4}[\# G_n]/\mathbf{E}_{4}[\# G_n \mid E_L]$ we can use  \eqref{first} to obtain $\mathbf{P}_{4}(E_L) \leq \kappa_{4}/n$ for some constant $\kappa_{4}>0$. By a symmetry argument, the same bound holds for the event $E_R=\{ \exists u \in G_n\,,\  \# (\op{D}(u) \cap \op{R}(u)) \geq c_3n/2\}$. Since $ \{G_{n}\ne \varnothing\}$ is the union of the events $E_R$, $E_L$ and $\{ \exists u \in G_n\,,\  \#\op{D}(u) \leq c_3n\}$, we easily deduce the upper bound of the theorem from the previous considerations and \eqref{troppetit}. \endproof

\section{Extensions} \label{extension} 
Fix $k\geq 2$. We can extend the recursive construction presented in the introduction by throwing polygons instead of chords: This will yield an analogue of the multitype branching process on the full $k$-ary tree. Formally if $x_{1}, \ldots , x_{k}$ are $k$ (distinct) points of $\mathbb{S}_{1}$ we denote by $\op{Pol}(x_{1}, \ldots,x_{k})$ the convex closure of $\{x_{1},\ldots,x_{k}\}$ in $\overline{\mathbb{D}}$. Let $(U_{i,j}: 1 \leq j \leq k \,,\ i \geq 1)$ be independent random variables that are uniformly distributed over $\mathbb{S}_1$. We construct inductively a sequence $L^{k}_1,L^{k}_2,\ldots$ of random closed subsets of the closed unit disk $\overline{\mathbb D}$.
To start with, $L^k_1$ is  $\op{Pol}(U_{1,1}, \ldots, U_{1,k})$. Then at step $n+1$, we consider two cases. Either the polygon $P_{n+1}:=\op{Pol}(U_{n+1,1}, \ldots ,U_{n+1,k})$ intersects $L^{k}_n$, and we put $L^k_{n+1}=L^k_n$. Or the polygon $P_{n+1}$ does not
intersect $L^k_n$, and we put $L^k_{n+1}=L^k_n \cup P_k$. Thus, for every integer $n\geq 1$, 
$L^{k}_n$ is a disjoint union of random $k$-gons. In a way very similar to what we did in the introduction we can identify the genealogy of the fragments appearing during this process with the complete $k$-ary tree
$$ \mathbb{T}_k \ =\  \bigcup_{i \geq 0} \{0,1, \ldots, k-1\}^i, \quad \mbox{ where } \{0,1, \ldots, k-1\}^0= \{\varnothing\}.$$ Then the number of ends of the fragments created during this process gives a labeling $\ell^k$ of $\mathbb{T}_k$ whose distribution can be described inductively by the following branching mechanism (this is an easy extension of \cite[Lemma 5.5]{CLG10}): For $u \in \mathbb{T}_{k}$ with label $m \geq 0$ we choose a decomposition  $m=m_1+m_2+ \ldots +m_k$ with $m_{1},m_{2}, \ldots, m_{k} \in \{0,1,\ldots,m\}$, uniformly at random among all  ${ m+k-1\choose k-1}$ possible choices, and we assign the labels $m_1+1, m_2+1, \ldots ,m_k+1$ to the children of $\varnothing$. Again the distribution of the labeling $\ell^{k}$ of $\mathbb{T}_k$ obtained  if we use the above branching mechanism but started from $a \geq 0$ at the root will be denoted by $\mathbf{P}_{a}$ and its expectation by $\mathbf{E}_{a}$. We use the same notation as in the binary case and are interested in a similar question: For which $a \geq 0$ does there exist with positive probability a ray $\mathbf{u}$ such that $\ell^k([\mathbf{u}]_i)  \geq a$ for every $i\geq 0$? Specifically, the value $a$ is called \emph{subcritical} for the process $(\ell^k(u),u \in \mathbb{T}_{k})$  when there exists a constant $c>0$ such that 
\begin{eqnarray*} \mathbf{P}_{a}( \exists u \in \{0,1,\ldots , k-1\}^n : \ell^k([u]_i) \geq a\,,\ \forall i \in \{0,1, \ldots , n\}) &\leq& \exp(-cn).\end{eqnarray*}
It is called \emph{supercritical} when there exists a constant $c>0$ such that we have both
$$\left\{\begin{array}{rcl} \mathbf{P}_{a}\left( \exists \mathbf{u} \in \{0,1,\ldots , k-1\}^\mathbb{N} : \ell^k([\mathbf{u}]_i) \geq a\,,\ \forall i \in \{0,1, \ldots \}\right) &\geq& c, \\
 \mathbf{E}_{a}[\# \left\{u \in \{0,1,\ldots , k-1\}^n : \ell^k([u]_i) \geq a\,,\ \forall i \in \{0,1, \ldots , n\}\right\}] &\geq& \exp(cn).\end{array}\right.$$
Note that a deterministic argument shows that if  $k\geq2$ and $a=2$, there always exists a ray with labels greater than or equal to $2$, also when $k=2$ and $a=3$ there exists a ray with labels greater than $3$. The case $k=2$ and $a=4$ has been treated in our main theorem. We have the following classification of all remaining cases: 
\begin{theorem} We have the following properties for the process $\ell^k$,
\begin{itemize} \item for $k=2$ and $a \geq 5$ the process is subcritical,
\item for $k=3$ the process is subcritical for $a\geq4$, and supercritical for $a=3$,
\item for $k \geq 4$ and $a \geq 3$ the process is subcritical.
\end{itemize}
\end{theorem}

\proof \textsc{Supercritical Case $k=3$ and $a=3$.} We will prove that for $k=3$ and $a=3$, the process  is supercritical. Similarly as in Section \ref{martingale} we consider the tree-indexed process $\ell^3$ on a fixed ray  of $\mathbb{T}_3$, say $\{0,0,0,\ldots \}$. Then the process $Y_n$ given by the $n$-th value of $\ell^3$ started from $3$ along this ray is a  homogeneous Markov chain with transition matrix given by 
\begin{eqnarray*} P_{3}(x,y) &=& \frac{2(x+2-y)}{(x+1)(x+2)}\mathbf{1}_{1\leq y \leq x+1}. \end{eqnarray*}
We introduce the stopping times $T_i = \inf\{n \geq 0, Y_n=i\}$ for $i=1,2$ and set $T= T_1\wedge T_2$. We consider a modification of the process $Y_n$ that we denote $\overline{Y}_n$, which has the same transition probabilities as $Y_n$ on $\{1,2,3,4\}$\,,\ but the transition between $4$ and $5$ for $Y_n$ is replaced by a transition from $4$ to $4$ for $\overline{Y}_n$. Thus we have $\overline{Y}_n \leq 4$ and an easy coupling argument shows that we can construct $Y_n$ and $\overline{Y}_n$ simultaneously in such a way that $\overline{Y}_n \leq Y_n$ for all $n \geq 0$. Hence we have the following stochastic inequality $$\overline{T} \ \leq \  T,$$ with an obvious notation for $\overline{T}$. To evaluate $\overline{T}$ we consider the subprocess $\overline{Y}_{n\wedge \overline{T}}$ which is again a Markov chain whose transition matrix restricted to $\{3,4\}$ is 
$$ \left( \begin{array}{cc}
{1/5} & {1/10}\\
{1/5} & {1/5}
\end{array} \right).$$
 The largest eigenvalue $\lambda_{\max}$ of this matrix is greater than $0.34$, which implies that $$ \mathbf{P}_{3}{T > n} \geq \mathbf{P}_{3}(\overline{T}>n) \geq \kappa_{5}(0.34)^n,$$ for some constant $\kappa_{5}>0$ independent of $n$. It follows that the expected number of paths starting at the root $\varnothing$ of $\mathbb{T}_3$ that have labels greater than or equal to $3$ up to level $n$, which is $3^n\P{T>n}$, eventually becomes strictly greater that $1$: There exists $n_{0} \geq 1$ such that $\mathbf{P}_{3}(T>n) > 3^{-n}$ for $n \geq n_{0}$. A simple coupling argument shows that the process $\ell^3$ started from $a \geq 3$ stochastically dominates the process $\ell^3$ started from $3$. Consequently, if we restrict our attention to the levels that are multiple of $n_{0}$ and declare that $v$ is a descendant of $u$ if along the geodesic between $u$ and $v$ the labels of $\ell^3$ are larger than $3$, then this restriction stochastically  dominates a supercritical Galton-Watson process. Hence the value $3$ is  supercritical for $\ell^3$.\\
 \textsc{Subcritical Case $k=3$ and $a=4$.}
As in the binary case we let 
$$ \tilde{P}_3(x,y) = \frac{2(x+2-y)}{(x+1)(x+2)}\mathbf{1}_{4\leq y \leq x+1},$$
be the substochastic matrix of the process $Y_n$ started at $4$ and killed when it hits $1,2$ or $3$. We will construct a positive vector $(h(x))_{x \geq 4}$ such that $\sum_x h(x) < \infty$ and 
\begin{eqnarray} h \cdot \tilde{P}_3 &\leq& \lambda h, \label{subcrit}\end{eqnarray} for some positive $\lambda < 1/3$, where we use the notation $h \cdot \tilde{P}_3(y) = \sum_x h(x) \tilde{P}_3(x,y)$. This will imply that \begin{eqnarray*}
\P{T>n} &\leq&  \frac{\sum_x h(x)}{h(4)}\lambda^n, \end{eqnarray*} where $T$ is the first hitting time of $\{1,2,3\}$ by the process $Y_n$ started at $4$. The subcriticality of the case $k=3$ and $a=4$ follows from the preceding bound since there are $3^n$ paths up to level $n$ and $\lambda <1/3$. To show the existence of a positive vector $x$ satisfying \eqref{subcrit} we begin by studying the largest eigenvalue of a finite approximation of the infinite matrix $\tilde{P}_3$. To be precise let $\tilde{P}_3^{(30)} = (\tilde{P}_3(i,j))_{4\leq i,j\leq 30}$. A numerical computation with Maple$^{\tiny \copyright}$ gives 
$$ \lambda_{\max}:=\max \left\{\mbox{Eigenvalues}(\tilde{P}_3^{(30)}) \right\} \simeq 0.248376642883065< 1/3.$$
The vector $(h(x))_{x \geq 4}$ is then constructed as follows. Let $(h(x))_{4 \leq x \leq 30}$ be an eigenvector associated with the largest eigenvalue $\lambda_{\max}$ of $\tilde{P}_{3}^{(30)}$, such that $\min_{4\leq x\leq 30}h(x)=h(30)=1$. Note that the vector $h$ can be chosen to have positive coordinates by the Perron-Frobenius theorem and it is easy to verify that $x \to h(x)$ is decreasing. For $x \geq 31$ we then let \begin{eqnarray*}
h(x) &=& 13^{x-30} \left(\frac{30!}{x!}\right)^{2}.\end{eqnarray*} We now verify that this vector satisfies \eqref{subcrit} with $\lambda$ slightly greater than $\lambda_{\max}$. Suppose first that $y \in \{4, \ldots, 30\}$. In this case $\sum_{4\leq x \leq 30} h(x) \tilde{P}_3(x,y)$ equals $\lambda_{\max} h(y)$ by definition, whereas the contribution of $\sum_{x \geq 31} h(x) \tilde{P}_3(x,y)$ is less than $\sum_{x\geq 31} h(x)<0.014$, thus 
\begin{eqnarray} h\cdot \tilde{P}_3(y) &\leq& 0.263 h(y). \end{eqnarray} Now, if $y \geq 31$ we have 
\begin{eqnarray*} \sum_{x \geq y-1} h(x)\tilde{P}_3(x,y) &\leq& 13^{y-30}\left(13^{-1}\tilde{P}_3(y-1,y) \left(\frac{30!}{(y-1)!}\right)^2  +  \tilde{P}_3(y,y)  \left(\frac{30!}{y!}\right)^2 
 +   \sum_{x \geq y+1} 13^{x-y}\left(\frac{30!}{x!}\right)^2 \right) \\
  &\leq& 13^{y-30} \left(\frac{30!}{y!}\right)^2 \left(\frac{2}{13}\frac{y^2}{y(y+1)} + \frac{4}{(y+1)(y+2)} + \sum_{i \geq 1} 13^i \left(\frac{y!}{(y+i)!}\right)^2\right) \\
   & \leq & 0.3\cdot h(y)
 \end{eqnarray*} which proves \eqref{subcrit}.\\
\textsc{Other critical cases.} The other critical cases are treated in the same way. We only provide the reader with the numerical values of the maximal eigenvalues of the truncated substochastic matrices that are very good approximations of the maximal eigenvalues of the infinite matrices,
$$\max \{\op{eigenvalues}(P_2(i,j))_{5\leq i,j\leq 30}\} \simeq 0.433040861268365 < 1/2,$$
$$\max \{\op{eigenvalues}(P_4(i,j))_{3\leq i,j\leq 30}\} \simeq 0.231280689028977 < 1/4.$$
\endproof

\bibliographystyle{abbrv}
\bibliography{/Users/nicolascurien/Dropbox/Macros-Bibli/bibli}

\end{document}